\documentclass[12pt]{amsart}
\usepackage{amscd}
\def\frk{\frak}               

\def\Phi{{\frk n}}
\def\Phi{{\frk N}}
\def\opn#1#2{\def#1{\operatorname{#2}}} 
\opn\chara{char}
\opn\length{\ell}
\opn\pd{pd}
\opn\rk{rk}
\opn\projdim{proj\,dim}
\opn\injdim{inj\,dim}
\opn\rank{rank}
\opn\depth{depth}
\opn\grade{grade}
\opn\height{height}
\opn\embdim{emb\,dim}
\opn\codim{codim}

\opn\Tr{Tr}
\opn\bigrank{big\,rank}
\opn\superheight{superheight}\opn\lcm{lcm}
\opn\trdeg{tr\,deg}%
\opn\reg{reg}
\opn\lreg{lreg}
\opn\ini{in}
\opn\div{div}
\opn\Div{Div}
\opn\cl{cl}
\opn\Cl{Cl}
\opn\Spec{Spec}
\opn\Supp{Supp}
\opn\supp{supp}
\opn\Sing{Sing}
\opn\Ass{Ass}
\opn\Ann{Ann}
\opn\Rad{Rad}
\opn\Soc{Soc}
\opn\Ker{Ker}
\opn\Coker{Coker}
\opn\Am{Am}
\opn\Hom{Hom}
\opn\Tor{Tor}
\opn\Ext{Ext}
\opn\End{End}
\opn\Aut{Aut}
\opn\id{id}

\opn\nat{nat}
\opn\pff{pf}
\opn\Pf{Pf}
\opn\GL{GL}
\opn\SL{SL}
\opn\mod{mod}
\opn\ord{ord}
\opn\aff{aff}
\opn\con{conv}
\opn\relint{relint}
\opn\st{st}
\opn\lk{lk}
\opn\cn{cn}
\opn\core{core}
\opn\vol{vol}
\opn\link{link}
\opn\star{star}
\opn\gr{gr}

\def\pot#1#2{#1[\kern-0.28ex[#2]\kern-0.28ex]}

\opn\dirlim{\underrightarrow{\lim}}
\opn\inivlim{\underleftarrow{\lim}}
\let\sect=\cap

\let\Dirsum=\bigoplus

\let\to=\rightarrow
\let\To=\longrightarrow
\def\Implies{\ifmmode\Longrightarrow \else
     \unskip${}\Longrightarrow{}$\ignorespaces\fi}
\def\implies{\ifmmode\Rightarrow \else
     \unskip${}\Rightarrow{}$\ignorespaces\fi}
\def\iff{\ifmmode\Longleftrightarrow \else
     \unskip${}\Longleftrightarrow{}$\ignorespaces\fi}

\let\:=\colon
\newtheorem{Theorem}{Theorem}[section]
\newtheorem{Lemma}[Theorem]{Lemma}
\newtheorem{Corollary}[Theorem]{Corollary}
\newtheorem{Proposition}[Theorem]{Proposition}

\let\epsilon\varepsilon
\let\phi=\varphi
\let\kappa=\varkappa
\def\qed{\ifhmode\textqed\fi
   \ifmmode\ifinner\quad\qedsymbol\else\dispqed\fi\fi}
\def\textqed{\unskip\nobreak\penalty50
    \hskip2em\hbox{}\nobreak\hfil\qedsymbol
    \parfillskip=0pt \finalhyphendemerits=0}
\def\dispqed{\rlap{\qquad\qedsymbol}}

\opn\dis{dis}
\def\pnt{{\raise0.5mm\hbox{\large\bf.}}}
\def\lpnt{{\hbox{\large\bf.}}}

\begin{document}

\title{Monomial ideals whose powers have a linear resolution}
\author{J\"urgen Herzog, Takayuki Hibi and Xinxian Zheng}
\address{J\"urgen Herzog, Fachbereich Mathematik und
Informatik,  
Universit\"at-GHS Essen, 45117 Essen, Germany}
\email{juergen.herzog@uni-essen.de}
\address{Takayuki Hibi, Department of Mathematics, Graduate School of Science,
Osaka University, Toyonaka, Osaka 560-0043, Japan}
\email{hibi@math.sci.osaka-u.ac.jp}
\address{Xinxian Zheng, Fachbereich Mathematik und
Informatik,  
Universit\"at-GHS Essen, 45117 Essen, Germany}
\email{xinxian.zheng@uni-essen.de}

\maketitle

\section*{Introduction}

In this paper we consider graded ideals in a polynomial ring over a field and ask when such an ideal has the property 
that all of its powers have a linear resolution.

It is known \cite{HT} that  polymatroidal ideals have linear resolutions and that powers of polymatroidal ideals are 
again polymatroidal (see \cite{CH} and \cite{HH}). In particular they have again linear resolutions. In general 
however, powers of ideals with linear resolution need not to have linear resolutions. The first example of such an 
ideal was given by Terai. He showed that over a base field of characteristic $\neq 2$ the Stanley Reisner ideal $I= 
(abd, abf, ace, adc, aef, bde, bcf, bce, cdf, def)$ of the minimal triangulation of the projective plane has a linear 
resolution, while $I^2$ has no linear resolution. The example depends on the characteristic of the base field. If the 
base field has characteristic 2, then $I$ itself has no linear resolution. 

Another example, namely $I=(def, cef, cdf, cde, bef, bcd, acf, ade)$ is given by Sturmfels \cite{St}. Again $I$ has a 
linear resolution, while $I^2$ has no linear resolution. The example of Sturmfels is interesting because of two 
reasons: 1.\ it does not depend on the characteristic of the base field, and 2.\ it is a linear quotient ideal. Recall 
that an equigenerated ideal $I$ is said to have {\em linear quotients} if there exists an order $f_1,\ldots, f_m$ of 
the  generators of $I$ such that for all $i=1,\ldots,m$ the colon ideals $(f_1,\ldots, f_{i-1}):f_i$ are generated by 
linear forms. It is quite easy to see that such an ideal has a linear resolution (independent on the characteristic of 
the base field). However the example of Sturmfels also shows  that powers of a linear quotient ideal need not to be 
again linear quotient ideals.  

On the other hand it is known (see \cite{CHT} and \cite{K}) that the regularity of powers $I^n$ of a graded ideal $I$ 
is bounded by a linear function $an+b$, and  is a linear function for large $n$. For ideals $I$ whose generators are 
all of degree $d$ one has the bound $\reg(I^n)\leq nd+\reg_x(R(I))$, as shown by R\"omer \cite{R}. Here $R(I)$ is the 
Rees ring of $I$ which is naturally bigraded, and $\reg_x(R(I))$  is the $x$-regularity of $R(I)$. It follows from 
this formula that each  power of $I$ has a linear resolution if $\reg_x(R(I))=0$. 

In this paper we will show (Theorem \ref{final}) that if $I\subset K[x_1,\ldots, x_n]$ is a monomial ideal with  
2-linear resolution, then each  power has a linear resolution.  Our proof is based on the formula of R\"omer. In the 
first section we give a new and very short proof of his result,  and remark that if there is a term order such that 
the initial ideal of the defining ideal $P$ of the Rees ring $R(I)$ is generated by monomials  which are linear in the 
variables $x_1,\ldots, x_n$, then $\reg_x(R(I))=0$.

In Section 2 we recall a result of Fr\"oberg \cite{F}  where he gives a combinatorial characterization of squarefree 
monomial ideals 2-linear resolution. A squarefree monomial ideal $I$ may be viewed as the edge ideal of a graph $G$.  
By Fr\"oberg, $I$ has a linear resolution if and only if the complementary graph $\overline{G}$ is chordal. There is 
an interesting characterization of chordal graphs due to G.A.Dirac \cite{D}. He showed that a graph is chordal if and 
only if it is the 1-skeleton of a quasi-tree. This characterization is essential for us in order to define the right 
lexicographical term order for which the initial ideal of $P$ is linear in the $x$ variables. We show this in the last 
section and use a description of the Graver basis of the egde ring of a graph due to Oshugi and Hibi 
\cite{OhHinormal}. Based on the same ideas and using polarization we also can treat monomial ideals which are not 
necessarily squarefree.

\section{The $x$-condition}
 
Let $K$ be a field, $S=K[x_1,\ldots, x_n]$ the polynomial ring, $I\subset S$ an equigenerated graded ideal, that is, a 
graded ideal whose generators $f_1,\ldots, f_m$ are all of same degree. Then the Rees ring
\[
R(I)=\Dirsum_{j\geq 0}I^jt^j=S[f_1t,\ldots, f_mt]\subset S[t]
\]
is naturally bigraded with $\deg(x_i)=(1,0)$ for $i=1,\ldots,n$ and $\deg(f_it)=(0,1)$ for $i=1,\ldots,m$. 

Let $T=S[y_1,\ldots,y_m]$ be the polynomial ring over $S$ in the variables $y_1,\ldots, y_m$. We define a bigrading on 
$T$ by setting $\deg(x_i)=(1,0)$ for $i=1,\ldots,n$, and $\deg(y_j)=(0,1)$ for $j=1,\ldots,m$. Then there is a natural 
surjective homomorphism of bigraded $K$-algebras $\phi\: T\to R(I)$ with $\phi(x_i)=x_i$ for $i=1,\ldots, n$ and 
$\phi(y_j)=f_jt$ for $j=1,\ldots,m$.

Let 
\[
\begin{CD}
\label{resolution}
F_\lpnt\: 0\to F_p@>>> F_{p-1}@>>> \cdots@>>> F_1@>>> F_0@>>> R(I)\to  0
\end{CD}
\]
be the bigraded minimal free $T$-resolution of $R(I)$. Here $F_i=\Dirsum_jT(-a_{ij},-b_{ij})$ for $i=0,\ldots, p$. The 
{\em $x$-regularity} of $R(I)$ is defined to be the number
\[
\reg_x(R(I))=\max_{i,j}\{a_{ij}-i\}.
\]

With the notation introduced one has the following result \cite[Theorem 5.3 (i)]{R} of R\"omer.

\begin{Theorem}
\label{R}
$
\reg(I^n)\leq nd+\reg_x(R(I)).
$
In particular, if $\reg_x(R(I))=0$, then each power of $I$ admits a linear resolution. 
\end{Theorem} 

For the reader's convenience we give a simple proof of this theorem: For all $n$, the exact sequence $F_\lpnt$ gives 
the exact sequence of graded $S$-modules
\begin{eqnarray}
\label{free}
\ \ \ \ \ 0\to (F_p)_{(*,n)}\To (F_{p-1})_{(*,n)}\cdots \To (F_1)_{(*,n)}\To (F_0)_{(*,n)}\To R(I)_{(*,n)}\to 0.
\end{eqnarray}
We note that $R(I)_{(*,n)}=I^n(-dn)$, and that $T(-a,-b)_{(*,n)}$ is isomorphic to the free $S$-module 
$\Dirsum_{|u|=n-b} S(-a)y^u$. It follows that (\ref{free}) is a (possibly non-minimal) graded free $S$-resolution of 
$I^n(-dn)$. This yields at once that $\reg(I^n(-dn))\leq \reg_x(R(I))$, and thus $\reg(I^n)\leq nd+\reg_x(R(I))$.

We say that $I$ satisfies the $x$-condition if $\reg_x(R(I))=0$.

\begin{Corollary}
\label{criterion}
Let $I\subset S$ be an equigenerated graded ideal, and let $R(I)=T/P$. Then each  power of $I$ has a linear resolution 
if for some term order $<$ on $T$ the defining ideal $P$ has a Gr\"obner basis $G$ whose elements  are at most linear 
in the variables $x_1,\ldots, x_n$, that is,  $\deg_x(f)\leq 1$ for all $f\in G$.    
\end{Corollary}

\begin{proof}
The hypothesis implies that $\ini(P)$ is generated by monomials $u_1,\ldots, u_m$ with $\deg_x(u_i)\leq 1$. Let 
$C_\lpnt$ be the Taylor resolution of $\ini(P)$. The module  $C_i$ has the  basis $e_\sigma$ with 
$\sigma=\{j_1<i_2<\ldots< j_i\}\subset [m]$. Each basis element $e_\sigma$ has the multidegree   $(a_\sigma,b_\sigma)$  
where $x^{a_\sigma}y^{b_\sigma}= \lcm\{u_{j_1},\ldots, u_{j_m}\}$. It follows that $\deg_x(e_\sigma)\leq i$ for all 
$e_\sigma\in C_i$. Since the shifts of $C_\lpnt$ bound the  shifts of a minimal multigraded  resolution of $\ini(P)$, 
we conclude that $\reg_x(T/\ini(P))=0$. On the other hand, by semi-continuity one always has $\reg_x(T/P)\leq 
\reg_x(T/\ini(P))$.
\end{proof}

\section{Monomial ideals with  2-linear resolution} 

Let $K$ be a field and $I\subset S=K[x_1,\ldots,x_n]$  be a squarefree monomial ideal generated in degree $2$. We may 
attach to $I$ a graph $G$ whose vertices  are the elements of  $[n]=\{1,\ldots,n\}$, and $\{i,j\}$ is an edge of $G$ 
if and only if $x_ix_j\in I$. The ideal   $I$ is called  the {\em  edge ideal} of $G$ and denoted $I(G)$. Thus the 
assignment $G\mapsto I(G)$ establishes a bijection between graphs and  squarefree monomial ideals generated in degree 
$2$.

The {\em complementary graph $\overline{G}$} of $G$ is the graph whose vertex set is again $[n]$ and whose edges are 
the non-edges of $G$. A graph $G$ is called {\em chordal} if each cycle of length $>3$ has a chord.  

We recall the following result of Fr\"oberg \cite[Theorem 1]{F} (see also \cite{V})

\begin{Theorem}[Fr\"oberg]
\label{froeberg} 
Let $G$ be graph. Then $I(G)$ has a linear resolution if and only if $\overline{G}$ is chordal.
\end{Theorem}

For our further considerations it is important to have a characterization of chordal graphs which is due to Dirac 
\cite{D}: let $\Delta$ be simplicial complex, and denote by ${\mathcal F}(\Delta)$ the set of facets of $\Delta$. A 
facet $F\in {\mathcal F}(\Delta)$  is called a {\em leaf} if either $F$ is the only facet of $\Delta$, or there exists  
$G\in {\mathcal F}(\Delta)$, $G\neq F$ such that $H\sect F\subset G\sect F$ for each $H\in{\mathcal F}(\Delta)$ with 
$H\neq F$. A vertex $i$  of $\Delta$ is called a {\em  free vertex} if $i$ belongs to precisely one facet. 

Faridi \cite{Fa} calls $\Delta$ a tree if each simplicial complex generated by a subset of the facets of $\Delta$ has 
leaf, and Zheng \cite{Z} calls $\Delta$ a {\em quasi-tree} if there exists a labeling $F_1,\ldots, F_m$ of the facets 
such that for all $i$ the facet $F_i$ is a leaf of the subcomplex $\langle F_1,\ldots, F_i\rangle$. We call such a 
labeling a {\em leaf order}. It is obvious that any tree is a quasi-tree, but the converse is not true. For us however 
the quasi-trees are important, because of 

\begin{Theorem}[Dirac]
\label{dirac}
A graph $G$ is chordal if and only if $G$ is the $1$-skeleton of a quasi-tree.
\end{Theorem}

As a consequence of Theorem \ref{froeberg} and \ref{dirac} we obtain

\begin{Proposition}
\label{order}
Let $I\subset S=K[x_1,\ldots,x_n]$ be a squarefree monomial ideal  with $2$-linear resolution. Then after suitable 
renumbering of the variables we have: if  $x_ix_j \in I$ with $i \neq j$, 
$k > i$ and $k > j$,
then either $x_ix_k$ or $x_jx_k$ belongs to $I$.
\end{Proposition}

\begin{proof}
We consider $I$ as the egde ideal of the graph $G$. Then by Theorem \ref{froeberg} and Theorem \ref{dirac} the 
complementary graph $\overline{G}$ is the $1$-skeleton of a quasi-tree $\Delta$. Let $F_1,\ldots, F_m$ be a leaf order 
of $\Delta$. Let $i_1$ be the number of free vertices of the leaf $F_m$. We label the free vertices of $F_m$ by 
$n,n-1,\ldots, n-i_1+1$, in any order. Next $F_{m-1}$ is a leaf of $\langle F_1,\ldots, F_{m-1}\rangle$. Say, 
$F_{m-1}$ has $i_2$ free vertices. Then we label the free vertices of $F_{m-1}$ by $n-i_1,\ldots, n-(i_1+i_2)+1$, in 
any order. Proceeding in this way we label all the vertices of $\Delta$, that is, those of $G$, and then choose the 
numbering of the variables of $S$ according to  this labeling. 

Suppose there exist $x_ix_j\in I$ and $k>i,j$ such that $x_ix_k\not\in I$ and
$x_jx_k\not\in I$. Let $r$ be the smallest number such that $\Gamma=\langle F_1,\ldots, F_r\rangle$ contains the 
vertices
$1,\ldots, k$. Then $k$ is a free vertex of $F_r$ in $\Gamma$. Since $x_ix_k\not\in I$ and 
$x_jx_k\not\in I$, we have that $\{i, k\}$ and  $\{j, k\}$ are edges of $\Gamma$, and since $k$ is a free vertex of 
$F_r$  in $\Gamma$ it follows that $i$ and $j$ are vertices of $F_r$.  Therefore  $\{i,j\}$ is an edge of $F_r$ and 
hence of $\Gamma$. However, this contradicts the assumption that  $x_ix_j\in I$. 
\end{proof}

We now consider a monomial ideal $I$ generated in degree $2$ which is not necessarily squarefree. Let $J\subset I$ be 
the ideal generated by all squarefree monomials in $I$. Then $I=(x_{i_1}^2,\ldots, x_{i_k}^2, J)$. 

\begin{Lemma}
\label{only} 
Suppose  $I$ has a linear resolution. Then $J$ has a linear resolution. 
\end{Lemma}

\begin{proof}
Polarizing (see \cite[Lemma 4.2.16]{BH}) the ideal $I=(x_{i_1}^2,\ldots, x_{i_k}^2, J)$ yields the ideal 
$I^*=(x_{i_1}y_1,\ldots, x_{i_k}y_{k}, J)$ in $K[x_1,\ldots,x_n, y_1,\ldots, y_k]$. We consider $I^*$ as the edge 
ideal of the graph $G^*$ with the vertices $-k,\ldots,-1, 1,\ldots,n$, where the vertices $-i$ correspond to the 
variables $y_i$ and the vertices $i$ to the variables $x_i$. Let $G$ be the restriction of $G^*$ to the vertices 
$1,\ldots, n$.  In other words, $\{i, j\}$ with $1\leq i<j$ is an edge of $G$ if and only it is an edge of $G^*$. Then 
it is clear that $J$ is the edge ideal of $G$. 

Assuming that $I$ has a linear resolution implies that $I^*$ has a linear resolution since $I^*$ is an unobstructed 
deformation of $I$. It follows that $\overline{G^*}$ is chordal, by Theorem \ref{froeberg}. Obviously the restriction 
of a chordal graph to a subset of the vertices is again chordal. Hence $\overline{G}$ is chordal, and so again by 
Theorem \ref{froeberg} we get that  $J$ has a linear resolution. 
\end{proof}

In the situation of  Lemma \ref{only} let $J=I(G)$, and let $\Delta$ be the quasi-tree whose  $1$-skeleton is 
$\overline {G}$, see Theorem \ref{froeberg} and Theorem \ref{dirac}.

\begin{Proposition}
\label{squares}
If $I=(x_{i_1}^2,\ldots, x_{i_k}^2, J)$ has a linear resolution, then  $i_j$ is a free vertex of $\Delta$ for 
$j=1,\ldots, k$, and no two of these vertices belong to the same facet. 
\end{Proposition}

\begin{proof}
We refer to the notation in the proof of Lemma \ref{only}.  Our assumption implies that $\overline{G^*}$ is chordal. 
Let $\Delta^*$ the quasi-tree whose $1$-skeleton is $\overline{G^*}$. 

Suppose that $i_j$ is not a free vertex of $\Delta$. Then there exist edges $\{i_j,r\}$ and $\{i_j,s\}$ in 
$\overline{G}$ such that $\{r,s\}$ is not an edge in $\overline{G}$. Then $\{i_j,r\}$ and $\{i_j,s\}$ are also edges 
in $\overline{G^*}$, and $\{r,s\}$ is not an edge in $\overline{G^*}$.  Since $x_{ij}y_j\in I^*$, it follows that 
$\{i_j,-j\}$ is not an edge in $G^*$, and since $x_ry_j$ and $x_sy_j$ do not belong to $I^*$ it follows that $\{-j, 
r\}$ and $\{-j, s\}$ are  edges of $\overline{G^*}$. Thus $\{i_j,r\}, \{r,-j\}, \{-j,s\}, \{s, i_j\}$ is circuit of 
length 4 with no chords, a contradiction. 

Suppose  $i_j$ and $i_l$  are free vertices belonging to the same facet of $\Delta$. Then $\{i_j, i_l\}$ is an edge in 
$\overline{G^*}$, and we also have that $\{i_j,-l\}$, $\{i_l,-j\}$ and $\{-j,-l\}$ are egdes of $\overline{G^*}$ since 
$x_{i_j}y_l$, $x_{i_l}y_j$ and $y_jy_l$ do not belong to $I^*$.  On the other hand, $\{i_j,-j\}$ and $\{i_l,-l\}$ are 
not edges of $\overline{G^*}$ since $x_{i_j}y_j$ and $x_{i_l}y_l$ belong to $I^*$. Therefore $\{i_j,i_l\}, 
\{i_l,-j\},\{-j,-l\},\{-l,i_j\}$ is the circuit of length 4 with  no chords, a contradiction.
\end{proof}

\begin{Corollary}
\label{true}
Suppose $I$ has a linear resolution and  $x_i^2\in I$. Then with the numbering of the variables as given in 
Proposition {\em \ref{order}} the following holds:  for all $j>i$ for which there exists $k$ such 
that $x_kx_j\in I$, one has $x_ix_j\in I$ or $x_ix_k\in I$.
\end{Corollary}

\begin{proof}
Suppose $x_i^2\in I$ and there exists a $j>i$ for which there exists $k$ such that $x_kx_j\in I$, but $x_ix_j$ and 
$x_ix_k$ both do not belong to  $I$. Then $k\neq i$, because $x_i^2\in I$. 

If $k\ne j$, then  $\{k,j\}$ is not an edge of $\Delta$, and $\{i,j\}$, $\{i,k\}$ both are
edges of $\Delta$. This implies that $i$ is not a free vertex of $\Delta$, contradicting  Proposition
\ref{squares}. 

If $k=j$, then $x_j^2\in I$ and $j$ is a free vertex of $\Delta$, by Proposition \ref{squares}. But since 
$x_ix_j\not\in I$ we have that  $\{i,j\}$ is an edge of $\Delta$. This implies that $i$ and $j$ belong to the same 
facet, again a contradiction to Proposition \ref{squares}.   
\end{proof}

\section{Monomial ideals satisfying the $x$-condition}

In the previous section we have seen that if $I$ is a monomial ideal generated in degree $2$ which has a linear 
resolution then it satisfies the conditions $(*)$ and $(**)$ listed in the next theorem. 

\begin{Theorem}
\label{takayuki}
Let $I \subset S = K[x_1, \ldots, x_n]$ be an ideal
which is generated by quadratic monomials and suppose that
$I$ possesses the following properties $(*)$ and $(**)$:

\medskip

$(*)$ if $x_ix_j \in I$ with $i \neq j$,
$k > i$ and $k > j$,
then either $x_ix_k$ or $x_jx_k$ belongs to $I$;

\medskip

$(**)$ if $x_i^2 \in I$ and $j > i$ for which there is
$k$ such that  $x_kx_j \in I$, then
either $x_ix_j \in I$ or $x_ix_k\in I$.

\medskip

\noindent
Let $R(I) = T/P$ be the Rees ring of $I$.
Then there exists a lexicographic order $<_{lex}$ on $T$
such that
the reduced Gr\"obner basis $G$ of the defining ideal $P$
with respect to $<_{lex}$ consists of binomials $f \in T$
with $\deg_{x}(f) \leq 1$.
\end{Theorem}

\begin{proof}
Let $\Omega$ denote the finite graph with the vertices
$1, \ldots, n , n + 1$ whose edge set $E(\Omega)$ consists of
those edges and loops
$\{ i, j \}$, $1 \leq i \leq j \leq n$, with
$x_i x_j \in I$ together with the edges
$\{ 1, n+1 \}, \{ 2, n+1 \}, \ldots, \{ n, n+1 \}$.
Let $K[\Omega] \subset S[x_{n+1}]$ denote
the {\em edge ring} of $\Omega$ studied in, e.g.,
\cite{OhHinormal} and \cite{OhHibinomial}.
%
%
%
%
%
%
Thus $K[\Omega]$ is the affine semigroup ring
generated by those quadratic monomials $x_i x_j$,
$1 \leq i \leq j \leq n + 1$,
with $\{ i, j \} \in E(\Omega)$.
Let $T = K[x_1, \ldots, x_n,
\{ y_{\{i, j\}} \}_{
1 \leq i \leq n, \, 1 \leq j \leq n
\atop \{i, j\} \in E(\Omega)}]$
be the polynomial ring and define the surjective
homomorphism
$\pi \, : \, T \to K[\Omega]$ by setting
$\pi(x_i) = x_i x_{n+1}$ and $\pi(y_{\{i, j\}}) = x_i x_j$.
The {\em toric ideal} of $K[\Omega]$
is the kernel of $\pi$.
Since the Rees ring $R(I)$ is isomorphic to
the edge ring $K[\Omega]$ in the obvious way,
we will identify
the defining ideal $P$ of the Rees ring with
the toric ideal of $K[\Omega]$.

We introduce the lexicographic order
$<_{lex}$ on $T$ induced by the ordering of the variables
as follows:
(i) $y_{\{i, j\}} > y_{\{p, q\}}$
if either
$\min\{i, j\} < \min\{p, q\}$
or
($\min\{i, j\} = \min\{p, q\}$
and $\max\{i, j\} < \max\{p, q\}$) and
(ii) $y_{\{i, j\}} > x_1 > x_2 > \cdots > x_n$
for all $y_{\{i, j\}}$.
Let $G$ denote the reduced Gr\"obner basis of
$P$ with respect to $<_{lex}$.

It follows (e.g., \cite[p. 516]{OhHibinomial}) that
the Graver basis of $P$ coincides with
the set of all binomials $f_{\Gamma}$,
where $\Gamma$ is a primitive even closed walk
in $\Omega$.
(In \cite{OhHibinomial}
a finite graph with no loop is mainly discussed.
However, all results obtained there are valid
for a finite graph allowing loops
with the obvious modification.)

Now, let $f$ be a binomial belonging to $G$ and
\[
\Gamma = (\{w_1, w_2\}, \{w_2, w_3\}, \ldots, \{w_{2m}, w_1\})
\]
the primitive even closed walk
in $\Omega$ associated with $f$.
In other words, with setting $y_{\{i,n+1\}} = x_i$
and $w_{2m+1} = w_1$,
one has
\[
f = f_{\Gamma} = \prod_{k=1}^{m} y_{\{w_{2k-1}, w_{2k}\}}-
\prod_{k=1}^{m} y_{\{w_{2k}, w_{2k+1}\}}.
\]
What we must prove is that,
among the vertices $w_1, w_2, \ldots, w_{2m}$,
the vertex $n+1$ appears at most one time.
Let $y_{\{w_{1}, w_{2}\}}$ be the biggest variable
appearing in $f$
with respect to $<_{lex}$ with $w_1 \leq w_2$.
Let $k_1, k_2, \ldots$ with $k_1 < k_2 < \cdots$ denote
the integers $3 \leq k < 2m$ for which $w_{k} = n + 1$.

\medskip

\noindent
{Case I:}
Let $k_1$ be even.
Since $\{n+1, w_1\} \in E(\Omega)$,
the closed walk
\[
\Gamma' = (\{w_1, w_2\}, \{w_2, w_3\}, \ldots,
\{w_{k_1-1}, w_{k_1}\},
\{w_{k_1}, w_1\})
\]
is an even closed walk in $\Omega$
with $\deg_x(f_{\Gamma'}) = 1$.  Since
the initial monomial $in_{<_{lex}}(f_{\Gamma'})
= y_{\{w_{1}, w_{2}\}} y_{\{w_{3}, w_{4}\}}
\cdots y_{\{w_{k_1-1}, w_{k_1}\}}$
of $f_{\Gamma'}$ divides
$in_{<_{lex}}(f_{\Gamma})
= \prod_{k=1}^{m} y_{\{w_{2k-1}, w_{2k}\}}$,
it follows that $f_{\Gamma} \not\in G$
unless $\Gamma' = \Gamma$.

\medskip

\noindent
{Case II:}
Let both $k_1$ and $k_2$ be odd.
This is impossible since $\Gamma$ is primitive
and since
the subwalk
\[
\Gamma'' = (
\{w_1, w_2\},
\ldots,
\{w_{k_1-1}, w_{k_1}\}, \{w_{k_2}, w_{k_2+1}\},
\ldots,
\{w_{2m}, w_{1}\}
)
\]
of $\Gamma$
is an even closed walk in $\Omega$.

\medskip

\noindent
{Case III:}
Let $k_1$ be odd and let $k_2$ be even.  Let $C$ be the
odd closed walk
\[
C = (\{w_{k_1}, w_{k_1+1}\},
\{w_{k_1+1}, w_{k_1+2}\},
\ldots,
\{w_{k_2-1}, w_{k_2}\})
\]
in $\Omega$.
Since both
$\{w_2, w_{k_1}\}$ and $\{w_{k_2}, w_{1}\}$
are edges of $\Omega$, the closed walk
\[
\Gamma''' = (\{w_1, w_2\}, \{w_2, w_{k_1}\},
C, \{w_{k_2}, w_{1}\})
\]
is an even closed walk in $\Omega$
and the initial monomial $in_{<_{lex}}(f_{\Gamma'''})$
of $f_{\Gamma'''}$ divides
$in_{<_{lex}}(f_{\Gamma})$.
Thus we discuss $\Gamma'''$ instead of $\Gamma$.

Since $\Gamma'''$ is primitive and since $C$ is of odd length,
it follows that
none of the vertices of $C$ coincides with
$w_1$ and that none of the vertices of $C$ coincides with
$w_2$.

%
%

(III -- a) First, we study the case when there is $p \geq 0$
with $k_1+p+2 < k_2$ such that $w_{k_1+p+1} \neq w_{k_1+p+2}$.
Let $W$ and $W'$ be the walks
\begin{eqnarray*}
W& = & (\{w_{k_1}, w_{k_1+1}\}, \{w_{k_1+1}, w_{k_1+2}\},
\ldots, \{w_{k_1+p+1}, w_{k_1+p+2}\}), \\
W' & = & (\{w_{k_2}, w_{k_2-1}\},
\{w_{k_2-1}, w_{k_2-2}\},
\ldots, \{w_{k_1+p+3}, w_{k_1+p+2}\})
\end{eqnarray*}
in $\Omega$.

(III -- a -- 1)
Let $w_1 \neq w_2$.
If either $\{w_2, w_{k_1+p+1}\}$ or
$\{w_2, w_{k_1+p+2}\}$
is an edge of $\Omega$, then
it is possible to construct an even closed walk
$\Gamma^{\sharp}$ in $\Omega$
such that
$in_{<_{lex}}(f_{\Gamma^{\sharp}})$ divides
$in_{<_{lex}}(f_{\Gamma'''})$
and $\deg_x(f_{\Gamma^{\sharp}}) = 1$.
For example, if, say,
$\{w_2, w_{k_1+p+2}\} \in E(\Omega)$
and if $p$ is even, then
\[
\Gamma^{\sharp} =
(\{w_2, w_1\}, \{w_1, w_{k_2}\}, W',
\{w_{k_1+p+2}, w_2\}
)
\]
is a desired even closed walk.

(III -- a -- 2)
Let $w_1 \neq w_2$.
Let
$\{w_2, w_{k_1+p+1}\} \not\in E(\Omega)$
and
$\{w_2, w_{k_1+p+2}\} \not\in E(\Omega)$.
Since $\{w_{k_1+p+1}, w_{k_1+p+2}\}$ is an edge of
$\Omega$, by $(*)$ either $w_2 < w_{k_1+p+1}$ or
$w_2 < w_{k_1+p+2}$.  Let $w_2 < w_{k_1+p+2}$.
Since $w_1 < w_2$ and $\{w_1, w_2\} \in E(\Omega)$,
again by $(*)$ one has
$\{w_1, w_{k_1+p+2}\} \in E(\Omega)$.
If $p$ is even, then consider the even closed walk
\[
\Gamma^{\flat} =
(
\{w_1, w_2\}, \{w_2, w_{k_2}\}, W',
\{w_{k_1+p+2}, w_1\}
)
\]
in $\Omega$.
If $p$ is odd, then consider the even closed walk
\[
\Gamma^{\flat} =
(
\{w_1, w_2\}, \{w_2, w_{k_1}\}, W,
\{w_{k_1+p+2}, w_1\}
)
\]
in $\Omega$.
In each case,
one has $\deg_x(f_{\Gamma^{\flat}}) = 1$.
Since $y_{\{w_1, w_2\}} > y_{\{w_1, w_{k_1+p+2}\}}$,
it follows that
$in_{<_{lex}}(f_{\Gamma^{\flat}})$ divides
$in_{<_{lex}}(f_{\Gamma'''})$.

(III -- a -- 3)
Let $w_1 = w_2$.
Since $w_1 < w_{k_1+p+1}$, by $(**)$
either
$\{w_1, w_{k_1+p+1}\} \in E(\Omega)$
or
$\{w_1, w_{k_1+p+2}\} \in E(\Omega)$.
Thus the same technique as in
(III -- a -- 2) can be applied.

(III -- b)
Second, if
$C = (\{n+1, j\}, \{j, j\}, \{j, n+1\})$,
then
in each of the cases
$w_1 < w_2 < j$,
$w_1 < j < w_2$ and
$w_1 = w_2 < j$,
by either $(*)$ or $(**)$,
one has either has
$\{w_1, j\} \in E(\Omega)$
or
$\{w_2, j\} \in E(\Omega)$.
\end{proof}

As the final conclusion of our considerations we obtain

\begin{Theorem}
\label{final}
Let $I$ be a monomial ideal generated in degree $2$. The following conditions are equivalent:
\begin{enumerate}
\item[(a)] $I$ has a linear resolution;
\item[(b)] $I$ has linear quotients;
\item[(c)] Each power of $I$ has a linear resolution.
\end{enumerate}
\end{Theorem}

\begin{proof} The implication (c)\implies (a) is trivial, while (b)\implies (a) is a general  fact. 
It follows from Proposition  \ref{order} and Corollary \ref{true} that if $I$ has a linear resolution, then the 
conditions $(*)$ and $(**)$ of Theorem \ref{takayuki} are satisfied, after a suitable renumbering of the variables. 
Hence by Corollary \ref{criterion} each power of $I$ has a linear resolution. 

It remains to prove (a)\implies (b): Again we may assume that the conditions $(*)$ and $(**)$ hold. Let $G(I)$ be the 
unique minimal set of monomial generators of $I$. We denote by $[u,v]$ the greatest common divisor of $u$ and $v$. 

We show that the following condition (q) is satisfied:  the elements of $G(I)$ can be ordered such  that if $u, v\in 
G(I)$ with $u>v$, then there exists $w>v$ such that $w/[w,v]$ is of degree $1$ and $w/[w,v]$ divides $u/[u,v]$.  This 
condition (q) then implies that $I$ has linear quotients.

The squarefree monomials in $G(I)$ will be ordered by the  lexicographical order induced by $x_n>x_{n-1}>\cdots >x_1$, 
and if $x_i^2\in G(I)$ then we let  $u>x_i^2>v$, where  $u$ is the smallest squarefree monomial of the form 
$x_kx_i$ with $k<i$, and where $v$ is the largest squarefree monomial less than $u$.  

Now, for any two monomials $u,v\in G(I)$ with $u>v$ corresponding to our order, we need to show that property $(q)$ 
holds. There are three cases:
 
Case 1: $u=x_sx_t$ and $v=x_ix_j$ both are squarefree monomials with $s<t$ and $i<j$. Since $u>v$, we have 
$t\geq j$. If $t=j$, take $w=u$. If $t>j$, then by $(*)$, either $x_ix_t\in G(I)$ or 
$x_jx_t\in G(I)$. Accordingly, let  $w=x_ix_t$ or $w=x_jx_t$. 

Case 2: $u=x_t^2$ and $v=x_ix_j$ with $i<j$. Since $u>v$, we have $t>j$. Hence by $(*)$, either 
$x_ix_t\in G(I)$ or $x_jx_t\in G(I)$. Accordingly, let $w=x_ix_t$ or $w=x_jx_t$. 

Case 3: $u=x_sx_t$ with $s\leq t$ and $v=x_i^2$. If $t=i$, then $s\neq t$ and take $w=u$. If $t>i$, then by 
$(**)$, we have either $x_ix_t\in G(I)$ or $x_ix_s\in G(I)$. Both elements  are greater than $v$ in our order.
Accordingly, let $w=x_ix_t$ or $w=x_ix_s$. Then again $(q)$ holds.

\end{proof}

\end{document}